\documentclass[12pt,reqno]{amsart}
\usepackage{enumerate, latexsym, amsmath, amsfonts, amssymb, amsthm, color}
\def\pmod #1{\ ({\rm{mod}}\ #1)}

\def\l{\left}
\def\r{\right}
\def\bg{\bigg}
\def\({\bg(}
\def\){\bg)}
\def\t{\text}
\def\f{\frac}

\def\ls{\leqslant}

\def\eq{\equiv}

\def\Proof{\noindent{\it Proof}}

\theoremstyle{plain}
\newtheorem{theorem}{Theorem}

\theoremstyle{definition}

\theoremstyle{remark}

\begin{document}
\title
[On Some determinant conjectures]
{On Some determinant conjectures}
\author[Ze-Hua Zhu and Chen-kai Ren]
{Ze-Hua Zhu and Chen-kai Ren}

\begin{abstract}
Let $p$ be a prime and  $c,d\in\mathbb{Z}$. Sun introduced the determinant $D_p^-(c,d):=\det[(i^2+cij+dj^2)^{p-2}]_{1<i,j<p-1}$ for $p>3$. In this paper, we confirm three conjectures on $D_p^-(c,d)$ proposed by Zhi-Wei Sun in \cite{1}.
\end{abstract}

\maketitle
\section{Introduction}
For an $n\times n$ matrix $[a_{ij}]_{1\leqslant i,j\leqslant n}$ over a commutative ring, we use $|a_{ij}|_{1\ls i,j\ls n}$ to denote its determinant.

Let $p$ be an odd prime, and let $(\frac{.}{p})$ be the Legendre symbol.
Carlitz \cite{C} determined the characteristic polynomial of the matrix
$$\left[x+\left(\frac{i-j}{p}\right)\right]_{1\leqslant i,j\leqslant p-1},
$$
and Chapman \cite{C1} evaluated the determinant
$$
\left| x+\left(\frac{i+j-1}{p}\right)\right|_{1\leqslant i,j\leqslant (p-1)/{2}}.
$$
Vsemirnov \cite{V12,V13} confirmed a challenging conjecture of Chapman
by evaluating the determinant
$$\left|\left(\frac{j-i}{p}\right)\right|_{1\leqslant i,j\leqslant (p+1)/{2}}.
$$
Sun \cite{S19} studied some determinants whose entries have the form $(\f{i^2+cij+dj^2}p)$,
where $c,d\in \mathbb{Z}$; in particular he introduced
$$(c,d)_p:=\left|\left(\frac{i^2+cij+dj^2}{p}\right)\right|_{1\ls i,j\ls p-1}$$
and $$[c,d]_p:=\left|\left(\frac{i^2+cij+dj^2}{p}\right)\right|_{0\ls i,j\ls p-1},$$
and proved that if $(\f dp)=1$ then
$$[c,d]_p=\begin{cases}\f{p-1}2(c,d)_p&\t{if}\ p\nmid c^2-4d,
\\\f{1-p}{p-2}(c,d)_p&\t{if}\ p\mid c^2-4d.\end{cases}$$
For any prime $p\eq3\pmod4$, Sun \cite[Remark 1.3]{S19} showed that
$$\left|\frac{1}{i^2+j^2}\right|_{1\leqslant i,j \leqslant (p-1)/2}\eq\l(\f 2p\r)\pmod p.
$$
For each prime $p\equiv 5\pmod{6}$, Sun \cite{S19} conjectured that $$2\left|\frac{1}{i^2-ij+j^2}\right|_{1\leqslant i,j \leqslant p-1}$$
 is a quadratic residue modulo $p$. This was recently confirmed by  Wu, She and Ni \cite{WSN}.

Let $p$ be an odd prime. For $b,c\in\mathbb{Z}$,
Sun \cite{S22} investigated the determinant
\begin{equation}D_p(b,c)=\left|(i^2+bij+cj^2)^{p-2}\right|_{1\leqslant i,j \leqslant p-1},
\end{equation}
and studied the Legendre symbol $(\f{D_p(b,c)}p)$.
Sun proved that if $p\eq2\pmod3,$ then
$D_p(1,1)=(\f{-2}p).$

Luo and Sun \cite{2} determined the Legendre symbol of $D_p(1,1)$ and $D_p(2,2)$ in some conditions.

She and Sun \cite{1} obtained some arithmetic properties of certain determinants involving
powers of $i^2+cij+dj^2$
, where $c$ and $d$ are integers.

Let $c,d\in\mathbb{Z}$. Sun \cite{1} introduced the determinant  $D_p^-(c,d):=\det[(i^2+cij+dj^2)^{p-2}]_{1<i,j<p-1}$ for any prime $p>3$. In this paper, we confirm some conjectures of Sun \cite[Conjecture 1.1]{1} as Theorem 1.1.
\begin{theorem}\label{1.1}
Let $p$ be a prime and $c,d\in\mathbb{Z}.$ The definition of $D_p^-(c,d)$ is as before. Then we have\\
(\romannumeral1)  $p\mid D_p^-(2,2)$ if $p\equiv 7\pmod 8$,\\
(\romannumeral2)  $p\mid D_p^-(3,3)$ if $p>5$ and $p\equiv 2\pmod 3$,\\
(\romannumeral3)  $p\mid D_p^-(3,1)$ if $p\equiv 3,7\pmod {20}$.
\end{theorem}
We will prove Theorem 1.1 in the next section.
\section{Proof of Theorem 1.1}
\setcounter{lemma}{0}
\setcounter{theorem}{0}
\setcounter{corollary}{0}
\setcounter{remark}{0}
\setcounter{equation}{0}
\Proof.  
(i) By Fermat's little theorem, there exists a unique polynomial
$$P(T)=a_0+a_1T+a_2T^2+\cdots+a_{p-2}T^{p-2}\in\mathbb{Z}[T]$$
such that 
$$(T^2+2T+2)^{p-2}\equiv P(T)\pmod p$$
for any $T\in\{1,2,\cdots,p-1\}$.
 
Now we recall the definition of the Lucas sequence. Given two integers $A$ and $B,$ the Lucas sequence $u_n=u_n(A,B)(n\in \mathbb{N})$ and its companion $v_n=v_n(A,B)(n\in \mathbb{N})$ are defined as follows:
$$u_0=0, u_1=1, u_{n+1}=Au_n-Bu_{n-1}(n=1,2,3,\cdots),$$
and
$$ v_0=2, v_1=A, v_{n+1}=Av_n-Bv_{n-1}(n=1,2,3,\cdots). $$
Let $u_k^{'}:=u_k(-2,2).$ By \cite[(4.11)]{2} we have
\begin{align}4a_k\equiv&(k+1)u^{'}_{p-k+1}-2(k-1)u^{'}_{p-k-1}\nonumber\\
&+2^{-k}((2k+4)u^{'}_k-ku^{'}_{k+2})\pmod p
\end{align}
where 
$$u_k^{'}=(-4)^{\lfloor\f k 4\rfloor}\times\begin{cases}
0 &\t{if}\ k\equiv 0\pmod 4,\\
1 &\t{if}\ k\equiv 1\pmod 4,\\
-2 &\t{if}\ k\equiv 2\pmod 4,\\
2 &\t{if}\ k\equiv 3\pmod 4.\end{cases}$$
Since $p\equiv 7\pmod 8$, it is easy to find that when $k\equiv 1\pmod 4$, $p-k-1\equiv 1\pmod 4$, $p-k+1\equiv 3\pmod 4$ and $k+2\equiv 3\pmod 4$. 
Hence, 
\begin{align}
4a_k\equiv&(k+1)\times2\times(-4)^{\lfloor\f{p-k+1}4\rfloor}-2(k-1)\times 1\times(-4)^{\lfloor\f{p-k-1}4\rfloor}\nonumber\\
&+2^{-k}\times[(2k+4)\times 1\times(-4)^{\lfloor\f k4\rfloor}-k\times 2\times(-4)^{\lfloor\f{k+2}4\rfloor}]\nonumber\\
\equiv&(k+1)\times2\times(-4)^{\f{p-k-2}4}-2(k-1)\times 1\times(-4)^{\f{p-k-2}4}\nonumber\\
&+2^{-k}\times[(2k+4)\times 1\times(-4)^{\f {k-1}4}-k\times 2\times(-4)^{\f{k-1}4}]\nonumber\\
\equiv &4\times(-4)^{\f{p-k-2}4}+2^{-k}\times4\times (-4)^{\f {k-1}4}\pmod p.
\end{align}

 Notice that $p-k-2\equiv k-1\pmod 8$ and $\left(\frac 2p\right)=1$. Then for $ k\equiv 1\pmod 4$ we have
\begin{align}
4a_k\equiv4\times(-1)^{\f{k-1}4}\times2^{\f{-k-1}2}\times(2^{\f{p-1}2}-1)\equiv&0\pmod p.
\end{align}
By \cite[Theorem 1.3]{1} and $(p-2)!\equiv 1\pmod p$,  we have
\begin{align}
D_p^-(2,2)=&\det[(i^2+2ij+2j^2)^{p-2}]_{1<i,j<p-1}\nonumber\\
=&\det[(j^{2(p-2)}( ij^{-1})^2+2(ij^{-1})+2)^{p-2}]_{1<i,j<p-1}\nonumber\\
\equiv&[(p-2)!]^{2(p-2)}\det[P(ij^{-1})]_{1<i,j<p-1}\nonumber\\
\equiv&4\sum_{i=0}^{\f {p-3}2}\hat{a_{2i}}\times\sum_{i=0}^{\f{p-3}2}\hat{a_{2i+1}}\pmod p,
\end{align}
where
\begin{align}
\hat{a_k}=\prod_{\begin{subarray} 0\le j\le p-2\\2|j-k, j\neq k\end{subarray}}a_j ,\qquad 
\end{align}
for  all $k=0,1,\cdots,p-2$.
By (2.2), we have $\hat{a_k}\equiv 0\pmod p$ when $2\nmid k$.
Then it follows that $D_p^-(2,2)\equiv 0 \pmod p$.

(ii) Similarly to (i), there exists a unique polynomial 
$$Q(T)=b_0+b_1T+b_2T^2+\cdots+b_{p-2}T^{p-2}\in\mathbb{Z}[T]$$
such that
$$(T^2+3T+3)^{p-2}\equiv Q(T)\pmod p$$
for any $T \in 1,2,\cdots,p-1$. 
 
Let $w_k:=u_k(-3,3).$ In view of Binet's formula, for any $k\in \mathbb{N}$ we have
$$w_k=\f{(\f {-3+\sqrt{3}i}{2})^k-(\f {-3-\sqrt{3}i}{2})^k}{\sqrt{3}i}$$
and thus
$$w_k=(-27)^{\lfloor\f k 6\rfloor}\times\begin{cases}
0 &\t{if}\ k\equiv 0\pmod 6,\\
1 &\t{if}\ k\equiv 1\pmod 6,\\
-3 &\t{if}\ k\equiv 2\pmod 6,\\
6 &\t{if}\ k\equiv 3\pmod 6,\\
-9&\t{if}\ k\equiv 4\pmod 6,\\
9&\t{if}\ k\equiv 5\pmod 6.\end{cases}$$
Similarly, by \cite[Lemma 2.1]{2} we have 
\begin{align}3b_k\equiv&(k+1)w_{p-k+1}-3(k-1)w_{p-k-1}\nonumber\\
&+3^{-k}((3k+6)w_k-kw_{k+2})\pmod p.
\end{align}
Since $p\equiv 2\pmod 3$, we noticed that when $k\equiv 2\pmod 6, p-k-1\equiv 2\pmod 6, p-k+1\equiv 4\pmod 6$ and $k+2\equiv 4\pmod 6$, then we have 
\begin{align}
    3b_k\equiv& -9\times(k+1)\times (-27)^{\f {p-k-3}6}+9\times(k-1)\times(-27)^{\f {p-k-3}6}\nonumber\\
    &+3^{-k}\times [-9\times(k+2)\times(-27)^{\f {k-2}6}+9\times k\times(-27)^{\f {k-2}6}]\nonumber\\
    \equiv&-18\times(-27)^{\f {p-k-3}6}+3^{-k}\times (-18)\times(-27)^{\f {k-2}6}\pmod p
\end{align}
Case 1: If $p\equiv 3\pmod 4$, then we have $p\equiv 11\pmod {12}$ and $\left(\frac 3p\right)=1$. Further more, we have $p-k-3\equiv k-2+6\pmod {12}$ when $k\equiv 2\pmod 6$, which means that $(-1)^{\f {p-k-3}6}=-(-1)^{\f {k-2}6}$. Then we have
\begin{align}
3b_k\equiv (-18)\times (-1)^{\f {k-2}6}\times (27)^{\f{-k-2}6}\times(1-3^{\f {p-1}6})\equiv 0\pmod p
\end{align}
Case 2: If $p\equiv 1\pmod 4$, then we have $p\equiv 5\pmod {12}$ and $\left(\frac 3p\right)=-1$. Further more, we have $p-k-3\equiv k-2\pmod {12}$ when $k\equiv 2\pmod 6$, which means that $(-1)^{\f {p-k-3}6}=(-1)^{\f {k-2}2}$. Then we have
\begin{align}
3b_k\equiv (-18)\times (-1)^{\f {k-2}6}\times (27)^{\f{-k-2}6}\times(1+3^{\f {p-1}2})\equiv 0\pmod p
\end{align}
In summary, we always have 
\begin{align}
3b_k\equiv 0\pmod p,
\end{align}
for all $k\equiv 2\pmod 6$.
Similarly to (2.4), we have
\begin{align}
D_p^-(3,3)\equiv4\sum_{i=0}^{\f {p-3}2}\hat{b_{2i}}\times\sum_{i=0}^{\f{p-3}2}\hat{b_{2i+1}}\pmod p,
\end{align}
where
\begin{align}
\hat{b_k}=\prod_{\begin{subarray} 0\le j\le p-2\\2|j-k, j\neq k\end{subarray}}b_j ,\qquad 
\end{align}
for all $k=0,1,\cdots,p-2$. By (2.10), we have we have $\hat{b_k}\equiv 0\pmod p$ when $2\mid k$.
Then it follows that $D_p^-(3,3)\equiv 0 \pmod p$.

(iii) Similarly to (i), there exists a unique polynomial 
$$R(T)=c_0+c_1T+c_2T^2+\cdots+c_{p-2}T^{p-2}\in\mathbb{Z}[T]$$
such that
$$(T^2+3T+1)^{p-2}\equiv R(T)\pmod p$$
for any $T \in 1,2,\cdots,p-1$. By \cite[Corollary 3.1]{S03}, $v_{\frac{p+1}{2}}(-3,1)=0$. By \cite{S22}, we find that if $v_{\frac{p+1}{2}}(-3,1)=0$, then 
$c_{\frac{p+1}{2}}=0.$ From \cite[Lemma 4.3]{2}, we obtain
\begin{equation}\label{ck}
   c_{\frac{p-3}{2}}=c_{p-1-\frac{p+1}{2}}= c_{\frac{p+1}{2}}=0.
\end{equation}
Notice that $\frac{p-3}{2}\equiv \frac{p+1}{2} \pmod{2}$. By  \cite[Theorem 1.3]{1}, we have $D_p^-(3,1)\equiv0 \pmod{p}.$

\end{document}